\documentclass[11pt]{amsart}
\usepackage{graphicx}
\usepackage[active]{srcltx}
 \makeatletter
\renewcommand*\subjclass[2][2000]{%
  \def\@subjclass{#2}%
  \@ifundefined{subjclassname@#1}{%
    \ClassWarning{\@classname}{Unknown edition (#1) of Mathematics
      Subject Classification; using '1991'.}%
  }{%
    \@xp\let\@xp\subjclassname\csname subjclassname@#1\endcsname
  }%
}
 \makeatother
\usepackage{enumerate,amssymb,  mathrsfs,yhmath}

\def\XXint#1#2#3{{\setbox0=\hbox{$#1{#2#3}{\int}$ }
\vcenter{\hbox{$#2#3$ }}\kern-.6\wd0}}

\newtheorem{theorem}{Theorem}[section]

\newtheorem*{lemma*}{Lemma}
\newtheorem{proposition}[theorem]{Proposition}
\newtheorem{corollary}[theorem]{Corollary}

\def\1ton{1,2,\ldots,n}

\def\ID{{\Bbb D}}

\usepackage{amssymb}
\usepackage{amsthm}
\usepackage{mathrsfs, amsfonts, amsmath}
\usepackage{graphicx}

\theoremstyle{definition}

\theoremstyle{remark}
\newtheorem{remark}[theorem]{Remark}

\numberwithin{equation}{section}

\def\XXint#1#2#3{{\setbox0=\hbox{$#1{#2#3}{\int}$}
\vcenter{\hbox{$#2#3$}}\kern-.5\wd0}}

\begin{document}

\title{Minimal surfaces and Schwarz lemma} \subjclass{Primary  53A10;
Secondary 42B30 }




\author{David Kalaj}
\address{University of Montenegro, Faculty of Natural Sciences and
Mathematics, Cetinjski put b.b. 81000 Podgorica, Montenegro}
\email{davidk@ucg.ac.me}

\begin{abstract}

We prove a sharp Schwarz type inequality for the  Weierstrass- Enneper representation of the minimal surfaces. It states the following. If $F:\mathbf{D}\to \Sigma$ is a conformal harmonic parameterization of a minimal disk $\Sigma$, where $\mathbf{D}$ is the unit disk and $|\Sigma|=\pi R^2$, then $|F_x(z)|(1-|z|^2)\le R$. If for some $z$ the previous inequality is equality, then the surface is an affine disk, and $F$ is linear up to a M\"obius transformation of the unit disk.

\end{abstract}  \maketitle


\section{Introduction}

The standard Schwarz-Pick lemma for holomorphic mappings states that every holomorphic mapping $f$ of the unit disk onto itself satisfies the inequality \begin{equation}\label{schar}
|f'(z)|\le \frac{1-|f(z)|^2}{1-|z|^2}.
\end{equation}
If the equality is attained in \eqref{schar} for a fixed $z=a\in\mathbf{D}$, then $f$ is a M\"obius transformation of the unit disk.

It follows from \eqref{schar} the weaker inequality  \begin{equation}\label{schar1}
|f'(z)|\le \frac{1}{1-|z|^2}
\end{equation}
with the equality in \eqref{schar1} for some fixed $z=a$ if and only if $f(z)=e^{it}\frac{z-a}{1-z\bar a}$. A certain extension of this result for harmonic mappings of the unit disk onto a Jordan domain has been given recently by the author in \cite{kalaj}.
We will extend this result to Weierstrass--Enneper parameterization
of minimal surfaces.

\subsection{Weierstrass--Enneper parameterization
of minimal surface}
The projections of minimal graphs in isothermal parameters are precisely
the harmonic mappings whose dilatations are squares of meromorphic functions.
 If $\Sigma$
is a minimal surface lying over a simply connected domain $D$ in
the $uv$ plane, expressed in isothermal parameters ($x$,
$y$), its projection onto
the base plane may be interpreted as a harmonic mapping $w =
f (z)$, where
$w =u +\imath v $ and $z =
x +\imath y.$ After suitable adjustment of parameters, it may be
assumed that $f$ is a sensepreserving harmonic mapping of the unit disk $\mathbf{D}$
onto, with $f (0) =w_0$ for some preassigned point $w_0$ in $D$. Let $f =
h +
\bar g$
be the canonical decomposition, where $h$ and $g$ are holomorphic. Then the dilatation
$\mu
=\frac{g'}{h'}$
of $f$ is an analytic
 function with  $|\mu(z)|<1$ in $\mathbf{D}$
and with the further property that $\mu =q^2$
for some function $q$ analytic in $U$. The minimal surface $\Sigma$
over $\Omega$
has the
isothermal representation $F=(u,v,t)$:

$$ u = \Re f (z) = \Re \int_0^z \phi_1(\zeta)d\zeta,$$

$$ v = \Im f (z) = \Im \int_0^z \phi_2(\zeta)d\zeta,$$

$$ t = \Im \int_0^z \phi_3(\zeta)d\zeta,$$

with $$\phi_1 =h'+ g'=p(1 + q^2),\ \phi_2 = -\imath (h'-
g') =-\imath p(1 -q^2), \ \ \ \text{and}\ \  \phi_3=
2ipq,$$
where $p$ and $q$ are the Weierstrass-Enneper parameters. Thus
$\phi_3^2 =-4\mu{h'}^2$
and
$h'= p$.
The first fundamental form of $\Sigma$
is $$ds^2 =\lambda^2 |dz|^2,$$ where
$$\lambda^2(z)=\frac{1}{2}\sum_{1}^{3}|\phi_k|^2.$$
A direct calculation shows that
$$\lambda=|h'|+ |g'|=|p|(1 + |q|^2).$$ For this fact and other important properties of minimal surfaces we refer to the book of Duren \cite{dure}.
Observe that $$|F_x|=|F_y|,\ \ \ \left<F_x,F_y\right>=0.$$
\section{The main results}
In this paper we consider minimal disks with Jordan boundaries and obtain some estimates of the conformal parametrization.

\begin{theorem}\label{mainth}

Let $F : \mathbf{D}\to \Sigma$
be the Weierstrass--Enneper parameterization
of Jordan minimal surface $\Sigma\subset
\mathbf{R}^3$ with the area $| \Sigma|=\pi R^2.$ Then the sharp
inequality

\begin{equation}\label{fx}
 |F_x(z)|\le \frac{R}{1-|z|^2}, \ \ \ z=x+\imath y\in\mathbf{D}
\end{equation}
holds. If for some $z$, the equality is attained, then $\Sigma$
is an affine disk and $F(m(x+\imath y))=F(0)+x N+y M$, where $M$ and $N$ are two orthogonal vectors of the equal length and $m$ is a M\"obius transformation of the unit disk onto itself. Moreover, every conformal mapping $F$ of the unit disk onto an affine disk of radius $R$ satisfies the equation in \eqref{fx} for some $a\in \mathbf{D}$.
\end{theorem}

 Thus  the previous theorem implies that
\begin{corollary}\label{rrje}
Let $F : \mathbf{D}\to \Sigma$
be the Weierstrass--Enneper parameterization
of Jordan minimal surface $\Sigma\subset
\mathbf{R}^3$ with the perimeter $| \partial \Sigma|=2\pi R.$ Then the sharp
inequality

\begin{equation}\label{fx2}
 |F_x|\le \frac{R}{1-|z|^2}, \ \ \ z=x+\imath y\in\mathbf{D}
\end{equation}
holds. If for some $z$, the equality is attained, then $\Sigma$
is an affine disk of radius $R$  and $F$ is composition of an affine mapping and a M\"obius transformation of the unit disk onto itself.

\end{corollary}

\begin{remark}
The same proof works for higher-dimensional case. So we can assume that the minimal disk is in $\mathbb{R}^n$ instead of $\mathbb{R}^3$.
\end{remark}
\begin{proof}[Proof of Corollary~\ref{rrje}]
Assume that $|\Sigma|=\mathrm{Area}(\Sigma)=\pi R_1^2$ and that $| \partial \Sigma|=2\pi R$. In view of isoperimetric inequality for minimal surfaces (see e.g. \cite{Osserman1978, top99, bol}  or a recent extension for harmonic surfaces \cite{kmm}) we have that  $$|\Sigma|\le \frac{| \partial \Sigma|^2}{4\pi}.$$ So $R_1\le R$. Further by \eqref{fx} we have
\begin{equation}\label{fx3}
 |F_x(z)|\le \frac{R_1}{1-|z|^2}, \ \ \ z=x+\imath y\in\mathbf{D}.
\end{equation}
So \eqref{fx2} follows at once.
\end{proof}
\begin{proof}[Proof of Theorem~\ref{mainth}]
We have $$F_x(z)=(\Re a'(z),\Re b'(z),\Re c'(z)).$$
Then $$|F_x(z)|^2=|\Re a'(z)|^2+|\Re b'(z)|^2+|\Re c'(z)|^2.$$ So $|F_x(z)|^2$ is subharmonic.

Thus by mean value inequality we have
$$|F_x(0)|^2\le \frac{1}{\pi}\int_{\mathbf{D}}|F_x|^2 dxdy=\frac{|\Sigma|}{\pi}.$$ Further let 
$$m(w)=\frac{w+z}{1+w\overline{z}}.$$ Then $$m'(0)=(1-|z|^2).$$ Define now theh mapping $H(z) = F(m(z))$. Then from the previous case we have \begin{equation}\label{hhh}|H_x(0)|\le \frac{|\Sigma|}{\pi}.\end{equation} Since $|H_x(0)|=|F_x(z)|(1-|z|^2)$, \eqref{hhh} implies that  $$|F_x(z)|^2\le \frac{1}{\pi}\frac{|\Sigma|}{(1-|z|^2)^2}=\frac{R^2}{(1-|z|^2)^2}.$$ This implies \eqref{fx}.

In order prove the equality statement, recall the definition of the Riesz measure $\mu$ of a subharmonic function $u$. Namely there exists a unique positive Borel measure $\mu$ so that $$\int_{\mathbf D} \varphi (z) d\mu(z)=\int_{\mathbf{D}} u \Delta \varphi(z) dm(z),\ \ \ \varphi\in C_0^2(\mathbf{D}).$$
Here $dm$ is the Lebesgue measure defined on the complex plane $\mathbf{C}$.
If $u\in C^2$, then $$\mathrm{d}\mu=\Delta u \mathrm{d}m.$$

\begin{proposition}\cite[Theorem 2.6 (Riesz representation theorem).]{pavlo}\label{pavlo}
If $u$ is a subharmonic function defined on the unit disk then for $r<1$ we have

\begin{equation}\label{subhar}\frac{1}{2\pi}\int_{\mathbf T} u(r z)|dz|-u(0)=\frac{1}{2\pi}\int_{|z|<r} \log\frac{r}{|z|} d\mu(z)\end{equation} where $\mu$ is the Riesz measure of $u$.
\end{proposition}
By applying Proposition~\ref{pavlo},  to the subharmonic function $$u(z)=|\Re a'(z)|^2+|\Re b'(z)|^2+|\Re c'(z)|^2$$ i.e. integrating \eqref{subhar} for $r\in[0,1]$, we obtain that \begin{equation}\label{kalfor}\frac{1}{2\pi}\int_0^1 r\int_{|z|<r} \log\frac{r}{|z|} d\mu(z)dr= \frac{1}{2\pi}\int_{\mathbf{D}} u(z)dm(z)-\frac{u(0)}{2}.\end{equation}
Assume first that the equality in \eqref{fx} is attained in $0$.
So if  $$|F_x(0)|^2= \frac{1}{\pi}\int_{\mathbf{D}}|F_x|^2 dxdy,$$ then the right-hand side of \eqref{kalfor} is zero.
Thus in particular we infer that $\mu=0$, or what is the same $\Delta u=0$. On the other hand $$\Delta u =|\nabla \Re a'(z)|^2+|\nabla \Re b'(z)|^2+|\nabla \Re c'(z)|^2.$$ So $\Re a'(z)$, $\Re b'(z)$ and $\Re c'(z)$ are constant functions, implying that \[\begin{split}F(z)&=(\Re a(z),\Re b(z), \Re c(z))\\&=F(0)+(a_0x+b_0y, c_0x+d_0y,e_0x+f_0y)=F(0)+x N+ y M\end{split}\] is an affine mapping. Here $$N=(a_0,c_0,e_0), \ \ \  M=(b_0, d_0, f_0).$$Further, conformality condition implies that $$R=|N|=|F_x(0)|=|F_y(0)|=|M|$$ and $$\left<N,M\right> =\left<F_x(0), F_y(0)\right>=0.$$ So $F(\mathbf{D})$ is a disk centered at $F(0)$ and with the radius $R$ lying in the plane $\mathcal{L}(N,M)$ spanned by the vectors $N$ and $M$.

If the equality statement is attained in $a\neq 0$, then consider the mapping $$H(z)= F\left(\frac{z+a}{1+z\bar{a}}\right).$$ Then $$|F'(a)| =\frac{R}{1-|a|^2},$$ if and only if $|H'(0)| =R$. Then by the previous proof $H$ is linear, and therefore $F(z) = H(\frac{z-a}{1-z\bar a})$.

Prove now the last part of the theorem. Assume now $\Sigma$ is an affine disk of radius $R$ ad let $F$ be a conformal mapping of $\mathbf{D}$ of $\Sigma$. Then there is a M\"obius transformation $m$ of the unit disk so that $$H(x+\imath y)=F\circ m(x+\imath y)=H(0)+(ax+by,cx+dy,ex+fy)=H(0)+x N+ y M$$ with $$R^2=|N|^2=a^2+c^2+e^2=b^2+d^2+f^2=|M|^2,\ \ \ \left<N, M\right>=ab+cd+ef=0.$$ Thus $$|H_x|=\sqrt{a^2+c^2+e^2}= R.$$ This implies that the equality in \eqref{fx} is attained for $H$ at $z=0$. Therefore $$|F_x(m(0))||m'(0)|=R.$$ Thus $$|F_x(a)|=\frac{R}{1-|a|^2}, \ \  a=m(0).$$

This concludes the proof.
\end{proof}

\begin{remark}
Recently in \cite{foka}, Forstneri\v c and the author have obtained, among the other result the following Schwarz lemma for conformal minimal parametrization of a minimal surface. Assume that $F:\mathbf{D}\to \mathbb{B}^n$ is conformal minimal parameterization of a minimal disk $\Sigma$, where $\mathbb{B}^n $ is the unit ball in $\mathbb{R}^n$. Assume also that $F(0)=0$. Then the following sharp inequality holds $|F(z)|\le |z|$ for $z\in\mathbf{D}$. If for some $z\neq 0$, $|F(z)|=|z|$, then $\Sigma=F(\mathbf{D})$ is an affine disk centered at $0$.
\end{remark}


\begin{thebibliography}{1}


\bibitem{bol}
G.~Bol.
\newblock Isoperimetrische {Ungleichungen} f{\"u}r {Bereiche} auf
  {Fl{\"a}chen}.
\newblock {\em Jahresber. Dtsch. Math.-Ver.}, 51:219--257, 1941.
\bibitem{dure}
P. Duren: {\it Harmonic mappings in the plane.} Cambridge University
Press, 2004.



\bibitem{foka} F. Forstneri\v c, D. Kalaj
\emph{Schwarz-Pick lemma for harmonic maps which are conformal at a point}, arXiv:2102.12403.


\bibitem{kmm}
D. Kalaj, M.  Markovic and M.  Mateljevic:  \emph{Carath\'eodory and Smirnov theorems for harmonic mappings of the unit disk onto surfaces}
Annales Academiae Scientiarum Fennicae
Mathematica
Volumen 38, 2013, 565--580.

\bibitem{kalaj}
D. Kalaj,  \emph{A sharp inequality for harmonic diffeomorphisms of the unit disk}, J. Geom. Anal. 29, No. 1, 392-401 (2019).
\bibitem{Osserman1978}
R.~Osserman.
\newblock The isoperimetric inequality.
\newblock {\em Bull. Amer. Math. Soc.}, 84(6):1182--1238, 1978.
\bibitem{pavlo} M. Pavlovi\'c: \emph{Function classes on the unit disc. An introduction.} De Gruyter
Studies in Mathematics, 52. De Gruyter, Berlin, 2014.


\bibitem{top99}
P.~Topping.
\newblock The isoperimetric inequality on a surface.
\newblock {\em Manuscr. Math.}, 100(1):23--33, 1999.
\end{thebibliography}
\end{document}